\documentclass[10pt]{article}

\usepackage{amsfonts,amsmath}
\usepackage{amsthm}
\usepackage{graphicx}
\usepackage{pdfsync}
\usepackage{epstopdf}

\usepackage{color}

 \newtheorem{thm}{Theorem}[section]
 \newtheorem{defn}[thm]{Definition}
 
 \newtheorem{cor}[thm]{Corollary}
 \newtheorem{lem}[thm]{Lemma}
 \newtheorem{prop}[thm]{Proposition}

\newtheorem{obs}[thm]{Observation}

\def\G{{\mathcal G}}

\def\field{{\mathbb F}}

\def\reals{{\mathbb R}}

\def\rmnmats #1 #2{M_{#1 #2}\,({\mathbb R})}
\def\cmnmats #1 #2{M_{#1 #2}\,({\mathbb C})}

\def\mnmats #1 #2{M_{#1 #2}}

\def\matdim #1 #2{#1 \times #2}

\def\rk{{\rm rank}}

\newcommand{\sS}{\mathcal{S}^-}
\newcommand{\smr}{\operatorname{mr}^-}

\newcommand{\mr}{\operatorname{mr}}

\newcommand{\x}{\times}
\newcommand{\bit}{\begin{itemize}}
\newcommand{\eit}{\end{itemize}}
\newcommand{\ben}{\begin{enumerate}}
\newcommand{\een}{\end{enumerate}}
\newcommand{\beq}{\begin{equation}}
\newcommand{\eeq}{\end{equation}}
\newcommand{\bea}{\begin{eqnarray*}}
\newcommand{\eea}{\end{eqnarray*}}
\newcommand{\bpf}{\begin{proof}}
\newcommand{\epf}{\end{proof}}

\title{A note on the minimum skew rank of powers of paths}

\author{Luz~M.~DeAlba\thanks{Department of Mathematics and Computer Science, Drake University, Des Moines, IA 50311, USA (luz.dealba@drake.edu).} \and Ethan~Kerzner\thanks{Department of Mathematics and Computer Science, Drake University, Des Moines, IA 50311, USA (ethan.kerzner@drake.edu).} \and Sarah~Tucker\thanks{Department of Mathematics and Computer Science, Drake University, Des Moines, IA 50311, USA (sarah.tucker@drake.edu).}.} 

\begin{document}

\maketitle

\begin{abstract}
The minimum skew rank of a simple graph $G$ is  the
smallest possible rank among all real skew symmetric  matrices, whose $(i,j)$-entry (for $i\ne j$) is nonzero whenever  $ij$ is an edge in $G$ and is zero otherwise.  In this paper we study the problem of minimum skew rank of powers and strict powers of paths.
\end{abstract}



\section{Introduction}\label{intro}

The minimum rank problem of a graph $G$ calls for the computation of the smallest possible rank among the matrices with a specific property (symmetric, skew-symmetric, Hermitian, positive definite), described by the graph $G$ and having entries in a given field $\field$. This problem has received considerable attention in the last few years (see for example~\cite{04BFH, 09BFHHHvdHS, 09DGHMR, 07FH, 09H}). Observe that many of these articles contain solutions to the minimum rank problem for special classes of graphs, matrices, or fields (see \cite{10IMA, 08BFS, 04BvdHL, 05BvdHL, 09BGL, 10DGKKMcDY, 10DS}). We have selected to work with skew-symmetric matrices over $\reals$, and the graphs under consideration are two different kinds of powers of paths.

The minimum rank of powers and strict powers of paths and trees was researched by the Minimum Rank Group at the AIM Workshop~\cite{AIM}, the results of this study appearing in~\cite{10DGKKMcDY}. The minimum skew rank problem was researched by the Minimum Skew Rank Group at the IMA Workshop~\cite{IMA}, with the first results appearing in~\cite{10IMA}, and further results in~\cite{10DS}. Throughout this paper, we adopt the notation and terminology from~\cite{10IMA},~\cite{10DGKKMcDY},~\cite{01W}, and~\cite{08F}.

A {\em graph} is a pair  $G = (V_G, E_G)$, where $V_G$ is the (finite, nonempty) set of vertices of $G$ and $E_G$ is the set of edges, where an edge is an unordered pair of vertices. All the graphs in this paper are {\em simple graphs} that is, all graphs are loop-free and undirected. The {\em order of a graph} $G$, denoted $|G|$, is the number of vertices of $G$. If  $e = uv \in E_G$, we say that $u$ and $v$ are {\em endpoints} of $e$; we also say that $u$ and $v$ are {\em adjacent}. Two graphs $G = (V, E)$ and $G' = \left( V', E' \right)$ are {\em isomorphic}, and we write $G \cong G'$,  whenever there exist bijections $\phi : V \rightarrow V'$ and $\psi : E \rightarrow E'$, such that $v \in V$ is an endpoint of $e \in E$ if and only if $\phi (v)$ is an endpoint of $\psi (e)$. A {\em subgraph} of a graph $G$ is a graph $H$ such that $V_H \subseteq V_G$ and $E_H \subseteq E_G$;  the graph $G - e$ denotes the subgraph $\left( V_G, E_G \setminus \{e\} \right)$ of $G$. If $W \subseteq V_G$ and $E' = \{ uv : u, v \in W, uv \in E_G \}$, the graph $(W, E')$ is referred to as the {\em subgraph of $G$ induced by $W$}. The subgraph of $G$ induced by $V_G \setminus \{v\}$ is denoted  by $G -v$. A {\em path} on $n$ vertices is the graph $P_n = \left( \{v_1,v_2, \dots, v_n \}, \{ e_i : e_i =  v_iv_{i+1}, 1 \le i \le n-1 \}   \right)$. A graph $G$, is {\em connected} if for every pair $u, v  \in V_G$, there is a path joining $u$ with $v$. A {\em walk of length $r$} in a graph $(V, E)$ is an alternating sequence: $v_{i_0}, e_{i_1}, v_{i_1}, e_{i_2},  \dots \ \dots, v_{i_r}, e_{i_r}$, of vertices, $v_{i_j} \in V$, and edges $e_{i_j} \in E$, (not necessarily distinct), such that $v_{i_{j-1}}$ and $v_{i_{j}}$ are the endpoints of $e_{i_j}$, for $j = 1, 2, \dots, r$. A {\em complete graph} is a graph whose vertices are pairwise adjacent, a complete graph on $n$ vertices is denoted by $K_n$. A graph $G$ is {\em bipartite} if $V_G = X \cup Y$, with $X \cap Y = \emptyset$, and such that each edge of $G$ has one endpoint in $X$ and the other in $Y$. A {\em complete bipartite graph}  is a bipartite graph in which each vertex in $X$ is adjacent to all the vertices in $Y$; a complete bipartite graph is denoted by $K_{n_1,n_2}$, where $|X| = n_1$ and $|Y| = n_2$. The {\em union} of graphs $G_1, G_2, \dots, G_k$, denoted $\bigcup_{i=1}^k \, G_i$, is the graph $\left( \cup_{i=1}^k \, V_i, \cup_{i=1}^k \, E_i \right)$.

A matrix $A \in \reals^{n \x n}$ is {\em skew-symmetric} if $A^T = -A$; note that the diagonal elements of a skew-symmetric matrix are zero. An $n \x n$ skew-symmetric matrix, $A = \left[ a_{ij} \right]$ is a {\em band} matrix of {\em bandwidth} $p$ if $a_{i,p+i} \ne 0, 1 \le i \le n-p$. The {\em graph of $A$}, denoted $\G (A)$, is the graph with vertices $\{1, . . . , n \}$ and edges $\{ ij : a_{ij} \ne 0, 1 \le i < j \le n \}$. Let $\sS (G) = \{ A \in \reals^{n \x n} : A^T = -A, \G(A) = G \}$ be the set of skew-symmetric matrices described by a graph  $G$. The {\em minimum skew rank of a graph $G$} is defined as $\smr (G) = \min \{ \rk(A) : A \in \sS (G) \}$.

\begin{defn}\cite[p.281]{00D}\label{upower}
Let $r$ be a positive integer and $G = \left( V_G, E_G \right)$ a graph. The graph $G$ to the power $r$ is the graph $G^r =  \left( V_G, E_{G^r} \right)$, where $ij \in E_{G^r}$ if and only if there is a walk in $G$ from vertex $i$ to vertex $j$ of length at most $r$. 
\end{defn}

\begin{defn}\label{spower}
Let $r$ be a positive integer and $G = \left( V_G, E_G \right)$ a graph. The graph $G$ to the strict power $r$ is the graph $G^{(r)} =  \left( V_G, E_{G^{(r)}} \right)$, where $ij \in E_{G^{(r)}}$ if and only if there is a walk in $G$ from vertex $i$ to vertex $j$ of length exactly $r$. 
\end{defn}

The following results are can be found in~\cite{10IMA}, where it was established that, in general minimum (symmetric) rank and minimum skew rank cannot be compared, in this paper we establish that if $G$ is a power of a path or a strict power of a path, then $\mr (G) \le \smr (G)$.

\begin{obs}\label{known}
Let $G$ be a graph.
\ben
\item If $H$ is an induced subgraph of $G$, then $\smr (H) \le \smr (G)$.
\item\label{components} If $G$ has connected components $G_1, \dots, G_k$, then $\smr (G) = \sum_{i=1}^k \, \smr \left( G_i \right)$.
\item\label{union} If $G = \bigcup_{i=1}^k \, G_i$, then $\smr (G) \le \sum_{i=1}^k \, \smr \left( G_i \right)$.
\een
\end{obs}

\begin{prop}\label{path_mr}\cite[Proposition 4.1]{10IMA}
For a path $P_n$ on $n$ vertices, 
$$\smr \left(P_n \right) = 
\left\{ \begin{array}{ll}
n & \mbox{if $n$ is even};\\
n-1 & \mbox{if $n$ is odd}.
\end{array} \right.$$
\end{prop}

%\begin{prop}\label{wheel}\cite[Proposition 3.11]{10IMA}
%If $F$ is an infinite field, $G'$ is connected, $|G| \ge 2$, and $G = G' \vee K_1$, then
%$\smr (F, G) = \smr(F, G')$.
%\end{prop}

\begin{thm}\label{mr2}\cite[Theorem 2.1]{10IMA} Let $G$ be a connected graph with $|G| \ge 2$ and let F be an infinite field. Then the following are equivalent:
\begin{enumerate}
\item $\smr (F, G) = 2$,
\item $G = K_{n_1,n_2,\dots, n_t}$ for some $t \ge 2$, $n_i \ge 1$, $i = 1, \dots, t$,
\item $G$ does not containt $P_4$ nor the paw as an induced subgraph.
\end{enumerate}
\end{thm}

\section{Minimum Skew Rank of Powers of Paths}\label{paths}

The connection between $\smr \left( P_n^r \right)$ and $\smr \left( P_n^{(r)} \right)$ is not made clear in~\cite{10DGKKMcDY}. However, we have found that when $r$ is even (see Theorem~\ref{reven}), these two quantities are closely related. For the remainder of our discussion we denote the set of vertices of $P_n$ by $\{ 1, 2, \dots, n-1, n \}$, where $1$ and $n$ are the pendant vertices, and $ij \in E_{P_n}$ if and only if $|i-j| = 1$. The following lemma contains some needed results, we omit the proof. 

\medskip

\begin{lem}\label{isomorphic1}
For a positive integer $m$, with $1 \le m \le n$, and $i \in \{ 1, 2, \dots, n-m+1 \}$, the induced subgraph of $P_n^r$ $\left( P_n^{(r)}, \mbox{ respectively} \right)$ on the set of vertices $\{ i, i+1, \dots, i+m-1\}$  is isomorphic to $P_m^r$ $\left( P_m^{(r)}, \mbox{ respectively} \right)$.
\end{lem}

\subsection{Minimum Skew Rank of Usual Powers of Paths}\label{usualp}

We know that $\smr \left(P_3 \right) = \smr \left(P_2 \right) = 2$, and that $P_2^r \cong K_2$, $P_3^r \cong K_3$, for $r \ge 2$, thus $\smr \left(P_3^r \right) = \smr \left(P_2^r \right) = 2$.

\begin{thm}\label{mainpathusual}
If $n$ and $r$ are positive integers, with $n \ge 4$, then
$$\smr \left(P_n^r \right) = 
\left\{ \begin{array}{llll}
n-r & \mbox{if } & 1 \le r \le n-3 \mbox{ and } n- r \mbox{ is even},\\
n-r + 1 & \mbox{if } & 1 \le r \le n-3 \mbox{ and } n-r \mbox{ is odd},\\
2 & \mbox{if} & r \ge n-2.
\end{array}
\right.$$
\end{thm}

\bpf
Recall (\cite{10DGKKMcDY}) that the graph $P_n^{n-2}$ is a complete multipartite graph isomorphic to  $K_n - e$, where $e$ is the edge $1n$. When $r \ge n-1$, the graph $P_n^{r}$ is a complete graph isomorphic to $K_n$. Thus if  $r \ge n-2$, it follows from Theorem~\ref{mr2}, that $\smr \left( P_n^{r} \right) = 2$.

Let $1 \le r \le n-3$, and observe that each matrix in $\sS \left(P_n^r \right)$ is the sum of $r+1$ band matrices of bandwidths $1, 2, \dots, r+1$ and zero diagonal,  thus the upper left $(n-r) \x (n-r)$ submatrix is lower-triangular with nonzero diagonal (\cite{08F}), and we have $\smr \left(P_n^r \right) \ge n - r$, except that, when $n-r$ is odd, we must have $\smr \left(P_n^r \right) \ge n - r+1$. In particular, note that if $r = n-3$ or $r = n-4$, then $\smr \left(P_n^r \right) \ge 4$.
 
\begin{figure}[h!]
\begin{center}
\scalebox{.3}{\includegraphics{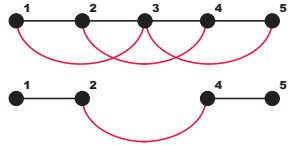}}
\caption{The graphs $P_5^2$ and $P_5^2 - \{ 3 \}$.}\label{wheelfig}
\end{center}
\end{figure}

For the cases $1 \le r \le n-3$, we proceed by induction on $n \ge 4$. For $n = 4$, we have  $\smr \left(P_4 \right) = 4 = (4 - 1) + 1$, and $\smr \left( P_4^2 \right) = 2$. For $n = 5$, we have  $\smr \left(P_5 \right) = 4 = 5 - 1$. Clearly (Figure~\ref{wheelfig}), $\smr \left(P_5^2 \right) = \smr \left( P_4 \right) = 4 = 5 - 2+1$, and form Theorem~\ref{mr2}$, \smr \left( P_5^3 \right) = 2$.

Let $k$ be a fixed integer, and suppose that for $4 \le n \le k-1$, 
$$\smr \left(P_{n}^r \right) = 
\left\{ \begin{array}{llll}
n - r & \mbox{if } &  1 \le r \le n - 3 \mbox{ and } n-r \mbox{ is even},\\
 n - r +1 &  \mbox{if } &  1 \le r \le n - 3 \mbox{ and } n-r \mbox{ is odd},\\
 2 & \mbox{if } &  r \ge n - 2. 
\end{array}
\right.$$

Let $n=k$, and let $r$ be an integer such that $1 \le r \le n-3$. Let $H_1$ be the subgraph of $P_n^r$, induced by the set of $r+2$ vertices $\{ 1, 2, \dots, r+2 \}$, and $H_2$  the subgraph of $P_n^r$,  induced by the set of $n-2$ vertices $\{ 3, 4, \dots, n \}$, so $H_1 \cong P_{r+2}^r$, and $H_2 \cong P_{n-2}^r$, and keeping the original labels on the vertices

$$P_n^r \cong \left( H_1 \cup \{r + 3, \dots, n-1, n \} \right) \bigcup \left( \{ 1, 2 \} \cup H_2 \right).$$

By~\ref{union}, in Observation~\ref{known}, 
$$\smr \left(P_n^r \right) \le \smr \left(H_1 \cup \{ r + 3, \dots, n-1, n \} \right) + \smr \left(  \{ 1, 2 \} \cup H_2 \right) =$$
$$\smr \left(H_1 \right) + \smr \left( H_2 \right) = \smr \left(P_{r+2}^r \right) + \smr \left( P_{n-2}^r \right).$$

By previous discussion, $\smr \left( P_{r+2}^r \right) = 2$, and by the induction hypothesis,
$$\smr \left( P_{n-2}^r \right) = 
\left\{ \begin{array}{lll}
n - r -2 & \mbox{if } & 1 \le r \le n -5 \mbox{ and } n- r -2 \mbox{ is even},\\
n- r - 1 & \mbox{if } & 1 \le r \le n-5 \mbox{ and } n-r -2 \mbox{ is odd},\\
2 & \mbox{if } & r \ge n-4.
\end{array}
\right.$$ 

Since $n$ and $n-2$ are both even or both odd, it follows that 
$$\smr \left(P_n^r \right) \le 2 + \left\{ \begin{array}{lll}
n - r -2 & \mbox{if } & 1 \le r \le n -5 \mbox{ and } n- r \mbox{ is even},\\
n- r - 1 & \mbox{if } & 1 \le r \le n-5 \mbox{ and } n-r \mbox{ is odd},\\
2 & \mbox{if } & r \ge n-4,
\end{array}
\right.$$
and consequently, that $$\smr \left(P_n^r \right) \le \left\{ \begin{array}{lll}
n - r & \mbox{if } & 1 \le r \le n -5 \mbox{ and } n- r \mbox{ is even},\\
n- r +1 & \mbox{if } & 1 \le r \le n-5 \mbox{ and } n-r \mbox{ is odd},\\
4 & \mbox{if } & r \ge n-4.
\end{array}
\right.$$
Since $\smr \left(P_n^r \right) \ge 4$ for $r = n-3$ or $r = n-4$, this is equivalent to 
$$\smr \left(P_n^r \right) \le \left\{ \begin{array}{lll}
n - r & \mbox{if } & 1 \le r \le n -3 \mbox{ and } n- r \mbox{ is even},\\
n- r +1 & \mbox{if } & 1 \le r \le n-3 \mbox{ and } n-r \mbox{ is odd}\\
4 & \mbox{if } & r \ge n-2
\end{array}
\right.$$
The proof of the theorem is now complete.
\epf

\begin{cor}
If $n$ and $r$ are positive integers, with $n \ge 3$, then $\mr \left( P_n^r \right) \le \smr \left( P_n^r \right)$, with equality if and only if $1 \le r \le n-3$ and $n - r$ is even. 
\end{cor}

\subsection{Minimum Skew Rank of Strict Powers of Paths}\label{strictp}

Note that $P_2^{(r)} \cong P_2 \cong K_2$, if $r$ is odd and $P_2^{(r)}  \cong K_1 \cup K_1$, if $r$ is even,  thus when $r \ge 1$, $\smr \left(P_2^{(r)} \right) = 2$ for $r$ odd, $\smr \left(P_2^{(r)} \right) = 0$ for $r$ even. Also, $P_3^{(r)} \cong P_3$ if $r$ is odd and $P_3^{(r)} \cong K_2 \cup K_1$ if $r$ is even, thus  for  $r \ge 1$, $\smr \left(P_3^{(r)} \right) = 2$.

%\newpage

\begin{lem}\label{strictterminal}
Let $n$ and $r$ be positive integers, with $n \ge 3$ and $r \ge n-2$. 
\ben 
\item If $r$ is odd, then $\smr \left(P_n^{(r)} \right) = 2$. 
\item If $r$ is even, then $\smr \left(P_n^{(r)} \right) = 4$.
\een
\end{lem}

\bpf
Recall (\cite{10DGKKMcDY}) that when $r \ge n-2$, and $r$ is odd,  $P_{n}^{(r)} \cong K_{\lfloor n/2 \rfloor,\lceil n/2 \rceil}$, and when $r$ is even,  $P_{n}^{(r)} \cong K_{\lfloor n/2 \rfloor} \cup K_{\lceil n/2 \rceil}$. The result follows from Theorem~\ref{mr2}.
\epf

\begin{thm}\label{rodd}
If $n$ and $r$ are positive integers, with $n \ge 3$ and $r$ odd, then
$$\smr \left(P_n^{(r)} \right) = 
\left\{ \begin{array}{lll}
n - r & \mbox{if } & 1 \le r \le n-3 \mbox{ and }  n \mbox{ is odd},\\
n- r +1 & \mbox{if } & 1 \le r \le n-3,  \mbox{ and }  n \mbox{ is even},\\
2 & \mbox{if } & r \ge n-2. 
\end{array}
\right.$$
\end{thm}

\bpf 
The case $r \ge n-2$ follows from Lemma~\ref{strictterminal}. Let $1 \le r \le n-3$, and observe that each matrix in $\sS \left(P_n^{(r)} \right)$ is the sum of $r+1$ band matrices (with zero diagonal) of bandwidths $1, 3, \dots, r+1$ if $r$ is even and bandwidths $2, 4, \dots, r+1$ if $r$ is odd,  thus the upper left $(n-r) \x (n-r)$ submatrix is lower-triangular with nonzero diagonal (\cite{08F}), thus the upper left $(n-r) \x (n-r)$ submatrix is lower-triangular with nonzero diagonal, and we have $\smr \left(P_n^{(r)} \right) \ge n - r$, except that, when $n-r$ is odd, we must have $\smr \left(P_n^{(r)} \right) \ge n - r+1$. In particular, note that if $r = n-3$ or $r = n-4$, then $\smr \left(P_n^{(r)} \right) \ge 4$.

Thus, $\smr \left(P_n^{(r)} \right) \ge
\left\{ \begin{array}{lll}
n - r & \mbox{if } & 1 \le r \le n-3 \mbox{ and }  n \mbox{ is odd},\\
n- r +1 & \mbox{if } & 1 \le r \le n-3 \mbox{ and } n \mbox{ is even}\\
2 & \mbox{if } & r \ge n - 2.
\end{array}
\right.$

The rest of the proof is by induction on $n$, note that the cases $n=3$ and $r \ge n-2$ follow from Lemma~\ref{strictterminal}. Let $1 \le r \le n-3$ and assume that for $3 \le n < k-1$ we have $\smr \left(P_n^{(r)} \right) = 
\left\{ \begin{array}{llll}
n - r & \mbox{if } & 1 \le r \le n-3, & n \mbox{ is odd},\\
n- r +1 & \mbox{if } & 1 \le r \le n-3, & n \mbox{ is even}\\
2 & \mbox{if } & r \ge n - 2.
\end{array}
\right.$ 

Let $n=k$, and let $H_1$ be the subgraph of $P_n^{(r)}$ induced by the set of vertices $\{ 1, 2, \dots, n-2 \}$, and $H_2$ the subgraph of $P_n^{(r)}$ induced by the set of vertices $ \{ n-r-1, n-r, \dots, n \}]$, so $H_1 \cong P_{n-2}^{(r)}$, and $H_2 \cong P_{r+2}^{(r)}$, and keeping the original labels on the vertices

$$P_n^{(r)} \cong \left( H_1 \cup \{ n-1, n \} \right) \bigcup \left( \{ 1, 2, \dots, n-r-2 \} \cup H_2 \right).$$

By~\ref{union}, in Observation~\ref{known},  
$$\smr \left(P_n^{(r)} \right) \le \smr \left( H_1 \cup \{ n-1, n \} \right) + \smr \left( \{ 1, 2, \dots, n-r-2 \} \cup H_2 \right) =$$ 
$$\smr \left( H_1 \right) + \smr \left( H_2 \right) = \smr \left( P_{n-2}^{(r)} \right) + \smr \left( P_{r+2}^{(r)} \right).$$
So by Lemma~\ref{strictterminal} and the induction hypothesis, $$\smr \left( P_{n-2}^{(r)} \right) = 
\left\{
\begin{array}{lll}
(n-2)-r & \mbox{if } & 1 \le r \le (n-2)-3  \mbox{ and }  n \mbox{ is odd},\\
(n-2) - r +1 & \mbox{if } & 1 \le r \le (n-2) -3 \mbox{ and }  n \mbox{ is even},\\
2 & \mbox{if } & r \ge (n-2) - 2
\end{array}
\right.$$
 and  $\smr \left( P_{r+2}^{(r)} \right) = 2$. It follows that 
$$\smr \left(P_n^{(r)} \right) \le 
\left\{
\begin{array}{llll}
(n-2 - r) + 2 & \mbox{if } & 1 \le r \le n-3 \mbox{ and }  n \mbox{ is odd},\\
(n-2 - r +1 ) +2 & \mbox{if } & 1 \le r \le n -3 \mbox{ and }  n \mbox{ is even},\\
4 & \mbox{if } & r \ge n-2
\end{array}
\right.$$ and hence, that $\smr \left(P_n^{(r)} \right) = 
\left\{
\begin{array}{llll}
n-r & \mbox{if } & 1 \le r \le n-3 \mbox{ and }  n \mbox{ is odd},\\
n - r +1 & \mbox{if } & 1 \le r \le n - 3 \mbox{ and }  n \mbox{ is even},\\
2 & \mbox{if } & r \ge n-2.
\end{array}
\right.$
\epf

\begin{lem}\label{isomorphism}
If $n$ and $m$ are positive integers with $n \ge 3$, then $P_n^{(2m)} \cong P_{\lfloor \frac{n}{2} \rfloor}^m \cup P_{\lceil \frac{n}{2} \rceil}^m $, the union being a disjoint union.
\end{lem}

\bpf
As before, we denote the set of vertices of $P_n$ by $\{ 1, 2, \dots, n-1, n \}$, the set of vertices of $P_{\lfloor \frac{n}{2} \rfloor}$ by $\{ v_1, v_2, \dots, v_{\lfloor \frac{n}{2} \rfloor} \}$, and the set of vertices of $P_{\lceil \frac{n}{2} \rceil}$ by $\{ u_1, u_2, \dots, u_{\lceil \frac{n}{2} \rceil} \}$. 

Define $\phi : \{ 1, 2, \dots, n-1, n \} \longrightarrow \{ v_1, v_2, \dots, v_{\lfloor \frac{n}{2} \rfloor} \} \cup \{ u_1, u_2, \dots, u_{\lceil \frac{n}{2} \rceil} \}$ by
$$\phi (w) = 
\left\{
\begin{array}{llll}
v_t & \mbox{if } & w = 2t, t \in \{ 1, 2, \dots, \lfloor \frac{n}{2} \rfloor \}\\
u_s & \mbox{if } & w=2s-1, s \in \{ 1, 2, \dots, \lceil \frac{n}{2} \rceil \}.
\end{array}
\right.$$
Clearly $\phi$ is an bijection. 

We know, from~\cite[Theorem~3.6]{10DGKKMcDY} (see also Theorem 8.1.3 in~\cite{08F}), that  $uv \in E_{P_n^{(2m)}}$ if and only if $|u-v| \in \{2m, 2m-2, 2m-4, \dots, 2 \}$, thus both $u$ and $v$ must be even, or both must be odd.
Define $\psi : E_{P_n^{(2m)}} \longrightarrow E_{P^m_{\lfloor \frac{n}{2} \rfloor}} \cup E_{P^m_{\lceil \frac{n}{2} \rceil}}$ by 
$$\psi \left( w_1 w_2 \right) = \left\{
\begin{array}{llll}
v_{t_1} v_{t_2} & \mbox{if } & w_1 = 2t_1,w_2 = 2t_2, t_1, t_2  \in \{ 1, 2, \dots, \lfloor \frac{n}{2} \rfloor \}\\
u_{s_1} u_{s_2} & \mbox{if } & w_1=2s_1-1, w_2=2s_2-1, s_1, s_2 \in \{ 1, 2, \dots, \lceil \frac{n}{2} \rceil \}.
\end{array}
\right.$$
Assume that $\psi \left( w_1 w_2 \right) = v_{t_1} v_{t_2}, $ for some $t_1, t_2  \in \{ 1, 2, \dots, \lfloor \frac{n}{2} \rfloor \}$, then $\left| v_{t_1} - v_{t_2} \right| = \left| t_1 - t_2 \right|  = \frac{1}{2} \left| w_1- w_2 \right| \in \{m, m-1, m-2, \dots, 1 \}$, thus $0 < \left| v_{t_1} - v_{t_2} \right| \le m$ and $v_{t_1} v_{t_2} \in E_{P^m_{\lfloor \frac{n}{2} \rfloor}}$. Similarly, $\psi \left( w_1 w_2 \right) = u_{s_1} u_{s_2} \in E_{P^m_{\lceil \frac{n}{2} \rceil}}$. 
The details for showing $\psi$ is an bijection are straightforward.
\epf

\begin{thm}\label{reven}
If $n$ and $r$ are positive integers, with $n \ge 4$, and $r=2s$, then
$$\smr \left(P_n^{(r)} \right) = 
\left\{ \begin{array}{lll}
n - r + 1& \mbox{if } & 1 \le r \le n-3 \mbox{ and }  n \mbox{ is odd},\\
n- r & \mbox{if } & 1 \le r \le n-3,  n = 2 t , t - s \mbox{ even}\\
n- r + 2 & \mbox{if } & 1 \le r \le n-3,  n = 2 t , t - s \mbox{ odd}\\
4 & \mbox{if } & r \ge n-2. 
\end{array}
\right.$$
\end{thm}

\bpf 
The case $r \ge n-2$ follows from Lemma~\ref{strictterminal}. Suppose $n$ is odd, and $1 \le r \le n-3$. By Lemma~\ref{isomorphism}, $P_n^{(r)} = P_n^{(2s)} \cong P_{\lceil \frac{n}{2} \rceil}^s \cup P_{\lfloor \frac{n}{2} \rfloor}^s$, and one of $\lceil \frac{n}{2} \rceil$ and $\lfloor \frac{n}{2} \rfloor$ is even and the other is odd. In either case, from Theorem~\ref{mainpathusual}, and~\ref{components} in Observation~\ref{known}, it follows that $\smr \left( P_n^{(r)} \right) = \smr \left( P_{\lceil \frac{n}{2} \rceil}^s \cup P_{\lfloor \frac{n}{2} \rfloor}^s \right) = \smr \left( P_{\lceil \frac{n}{2} \rceil}^s  \right) + \smr \left( P_{\lfloor \frac{n}{2} \rfloor}^s \right) = \lceil \frac{n}{2} \rceil - s + \lfloor \frac{n}{2} \rfloor - s + 1 = n -2s + 1 = n - r +1$.

Now let $n$ be even, $n = 2t$, and $1 \le r \le n-3$, so that by Lemma~\ref{isomorphism}, $P_n^{(r)} = P_n^{(2s)} \cong P_t^s \cup P_t^s$, and thus $\smr \left( P_n^{(r)} \right) = 2 \smr \left( P_t^s \right)$. From Theorem~\ref{mainpathusual}, $\smr \left( P_t^s \right) = t - s$, if $t - s$ is even, and $t - s +1$, if $t - s$ is odd. Therefore, $\smr \left( P_n^{(r)} \right) = n-r$, if $t - s$ is even, and $n - r +2$, if $t - s$ is odd.
\epf

\begin{cor}
If $n$ and $r$ are positive integers with $n \ge 3$, then $\mr \left( P_n^{(r)} \right) \le \smr \left( P_n^{(r)} \right)$, with equality if and only if $n$ and $r$ are both odd, or $n$ and $r$ are both even, $n = 2t$, $r = 2s$ and $t-s$ is even.
\end{cor}
%
%{\bf Acknowledgement.} Luz M. DeAlba began research on minimum rank of powers of graphs while at the  American Institute of Mathematics~\cite{AIM}. She began research on minimum skew rank under the sponsorship of the Institute of Applied Mathematics and Iowa State University~\cite{IMA}.

\end{document}